\documentclass[leqno, 11pt]{article}
\usepackage[utf8]{inputenc}

\usepackage{amsmath}
\usepackage{cleveref}

\usepackage{lipsum}
\usepackage{amsfonts}
\usepackage{amssymb}
\usepackage{graphicx}
\usepackage{epstopdf}
\usepackage{algorithm, algorithmic}
\usepackage{verbatim}
\ifpdf
  \DeclareGraphicsExtensions{.eps,.pdf,.png,.jpg}
\else
  \DeclareGraphicsExtensions{.eps}
\fi

\graphicspath{ {./images/} }
\usepackage [english]{babel}
\usepackage [autostyle, english = american]{csquotes}
\MakeOuterQuote{"}

\title{Hamiltonian Monte Carlo methods for spectroscopy data analysis}

\author{Daniel McBride \thanks{Department of Mathematics, University of Tennessee Knoxville 
  (dmcbride@utk.edu).}
\and Ioannis Sgouralis\thanks{Department of Mathematics, University of Tennessee Knoxville 
  (isgoural@utk.edu).}
  }

\newcommand{\HN}{\mathcal{N}}

\newcommand{\HT}{\mathcal{T}}

\newcommand{\R}{\mathbb{R}}

\def\l{\lambda}

\newcommand{\Ker}[1]{\mathsf{Ker}~}

\def\l{\lambda}

\def\m{\mu}

\def\and{\quad\text{and}\quad}

\def\t{\tau}
\def\o{\omega}

\def\Iref{I_{\text{ref}}}

\def\Iref{I_{\text{ref}}}
\def\Ibg{I_{\text{bg}}}
\def\Vlike{V_{\text{like}}}
\def\Vpri{V_{\text{prior}}}

\def\Hpri{H_{\text{prior}}}
\def\Hlike{H_{\text{like}}}
\def\Tpri{T_{\text{prior}}}
\def\Tlike{T_{\text{like}}}
\def\grad{\nabla}

\def\o{\omega}
\def\texp{\tau_\text{exp}}
\def\tdead{\tau_\text{dead}}
\def\tsub{\tau_\text{sub}}

\def\th{\theta}

\usepackage[margin=1in,footskip=0.25in]{geometry}

\date{}

\begin{document}

\maketitle

\begin{abstract}
We present a scalable Bayesian framework for the analysis of confocal fluorescence spectroscopy data, addressing key limitations in traditional fluorescence correlation spectroscopy methods. Our framework captures molecular motion, microscope optics, and photon detection with high fidelity, enabling statistical inference of molecule trajectories from raw photon count data, introducing a superresolution parameter which further enhances trajectory estimation beyond the native time resolution of data acquisition. To handle the high dimensionality of the arising posterior distribution, we develop a family of Hamiltonian Monte Carlo (HMC) algorithms that leverages the unique characteristics inherent to spectroscopy data analysis. Here, due to the highly-coupled correlation structure of the target posterior distribution, HMC requires the numerical solution of a stiff ordinary differential equation containing a two-scale discrete Laplacian. By considering the spectral properties of this operator, we produce a CFL-type integrator stability condition for the standard St{\"o}rmer-Verlet integrator used in HMC. To circumvent this instability we introduce a semi-implicit (IMEX) method which treats the stiff and non-stiff parts differently, while leveraging the sparse structure of the discrete Laplacian for computational efficiency. Detailed numerical experiments demonstrate that this method improves upon fully explicit approaches, allowing larger HMC step sizes and maintaining second-order accuracy in position and energy. Our framework provides a foundation for extensions to more complex models such as surface constrained molecular motion or motion with multiple diffusion modes.

\end{abstract}

\


\


\section{Introduction}
Fluorescence spectroscopy studies light-matter interactions by means of specialized experiments \cite{B2006, L2006, S2023I, S2023II, S2023III}. Nowadays, various spectroscopic modalities are routinely used to identify chemical compounds \cite{K2020, M1997}, characterize materials \cite{Y2011}, image biological specimens \cite{L2007,L2015}, and, most recently, also to quantify dynamics of moving particles \cite{J2019, J2019method, S2023, S2017}. In particular, \emph{fluorescence correlation spectroscopy} (FCS) combines unique optical devices and sampling electronics to record direct data at time scales on the order of one nanosecond which can resolve single photon emissions \cite{L2006, T2020, S2023I, S2023II, S2023III}. Such raw data, combined with mathematical analysis \cite{J2019, E1974, M1974, D2020}, elucidate dynamics of moving molecules in chemical,  biological and biochemical applications.

A FCS experiment uses a confocal microscope to illuminate only a diffraction limited spot within the specimen of interest \cite{W1996,L2007,ZZOM2007}. Subsequently, fluorescent particles, that absorb the illuminating light and re-emit it at a different color, are allowed to move within this volume. Optical filters are used to separate the emitted light and fast photon detectors are then employed to record them in real time. This way, a FCS experiment yields time series data of photon detections that are acquired typically at rates as fast as one sample per nano- to millisecond for time periods in excess of one minute. The resulting time series datasets are subsequently analyzed to estimate mainly diffusion coefficients that are characteristic of the motion of the fluorescent particles \cite{J2019}.

Data acquisition rates at the nanosecond scale across typical durations at or longer than a minute produce datasets with sizes on or above the order of \( 10^{10} \). Processing such high data volume requires careful design and implementation of statistical analysis techniques.

Traditionally, data analysis in FCS is carried out by means of correlation functions that are derived based on simplistic theoretical models and fitted heuristically. Among many FCS simplifications the most important ones ignore detection noise and optics \cite{E1974,M1974,R1993}. As shown in \cite{J2019}, such an approach has a degrading effect on the quality of the generated estimates that may be misleading or biased.

In contrast, direct analysis methods using Markov state-space models \cite{BN2006,P2023} in combination with Bayesian non-parametrics \cite{G2017, P2023} have recently been proposed as viable alternatives to correlative FCS analysis \cite{J2019,D2020}. By relying on realistic modeling representations without artificial simplifications, direct approaches yield reliable estimates with only a fraction of the data traditionally required for FCS applications. This is of paramount importance for biological FCS applications which seek to reduce photo-damage and so cannot afford excessively long data acquisition periods or exceedingly high illumination intensities as required for high signal-to-noise ratio (SNR) \cite{B1991}. Furthermore, such specialized approaches are built on high fidelity representations of the experiment and allow for modifications and generalizations to non-standard setups. For instance, a faithful representation allows for modeling read-out noise and detectors \cite{K2021extraction, K2021}, complicated optics \cite{G2005, L2007, S2019}, multiple confocal volumes \cite{J2022}, as well as photo-physics \cite{K2021, S2018} which are common limiting factors to conventional FCS.

The study in \cite{J2019} provided a proof of principle for the direct analysis of raw spectroscopic data acquired in a typical FCS experiment. However, the numerical methods in \cite{J2019} are limited to slowly diffusing particles ($\lesssim$10~$\mu$m$^2$/sec) that may not be always given. In fact, biomolecules assessed in vitro or data of low SNR such as free fluorescent dyes and inorganic compounds assessed in vivo may exhibit diffusion coefficients as high as $\approx$500~$\mu$m$^2$/sec \cite{S2017,K2021,S2017}. Extracting estimates of such high diffusion coefficients from raw spectroscopic data requires a more powerful computational framework than that proposed in \cite{J2019}.

In this study, we develop and investigate a novel framework for analyzing confocal fluorescence spectroscopy data. Specifically, we construct a generative probabilistic model of this spectroscopy method and examine the numerical aspects of Bayesian techniques used to characterize the conditional distribution of molecule trajectories given integrative photon count data \cite{K2021,K2021extraction}. This requires solving a high-dimensional Bayesian inverse problem involving a target distribution exhibiting a complex covariance structure. Given the dimensionality of this distribution (detailed in \cref{sec:methods}), standard Monte Carlo approaches, such as Random Walk Metropolis-Hastings, are unsuitable \cite{L2001, RC1999}. Instead, we put forward a Hamiltonian Monte Carlo (HMC) scheme \cite{B2017, BRSS2018, GJM2011}, which efficiently handles high-dimensional distributions. However, the computational demands are high because our approach requires repeatedly solving an ordinary differential equation (ODE). For typical fluorescence spectroscopy applications, the ODE is stiff, necessitating the use of ill-conditioned matrix operations for accurate numerical solutions. Our problem, typical of Bayesian inversion problems, involves a posterior distribution formed as the product of a likelihood function and a prior \cite{P2023,BN2006,RC1999,L2001}. This multiplicative structure naturally separates the ODE into stiff and non-stiff components. We analyze a solver that leverages this split structure and propose a new HMC-based inference method for fluorescence spectroscopy data, validated using in silico data mimicking realistic confocal fluorescence spectroscopy conditions.

The remainder of this study is organized as follows. In \cref{sec:methods} we describe our model, the forward problem of generating synthetic data, and the inverse problem of estimating the diffusion coefficient of a molecular species given photon count data. The motivation for and description of our new HMC algorithms are also in \cref{sec:methods}. We provide demonstrations and numerical
results in \cref{sec:results}, and a discussion in
\cref{sec:discussion}.

\section{Methods}
\label{sec:methods}

In this section, we first describe a mathematical model of fluorescence spectroscopy. For clarity of exposition, we focus on a motion model in one spatial dimension for homogeneous single molecule diffusion. Subsequently, we use this model to develop a statistical method that estimates molecular trajectories from raw spectroscopic measurements. Finally, we use the same model to investigate the numerical properties of our algorithms.

\subsection{Model description}
\label{sec:model}

We represent each molecule as a one-dimensional material point and denote its position by \( x\in\R \). The photon emission rate of a molecule at position \( x \) is governed by
\begin{align}
    I(x)=\Ibg+\Iref G(x).
\end{align} 
The quantity \( \Ibg>0 \) is the background photon emission rate which we assume remains constant thorough the experiment. Further, $\Iref>0$ determines the reference photon rate, also known as the intensity or brightness \cite{L2006}, of the \emph{point spread function} \( G(x) \). 

The point spread function $G(x)>0$ is unitless and describes the distribution of light emitted by a molecule at position $x$ in the confocal volume due to the optical setup \cite{G2005, L2007, S2019, W1996, ZZOM2007}. Following existing approaches, we model the point spread function by a Gaussian centered at the origin,
\begin{align}
    G(x)=\exp\left(-\frac{x^2}{2\o}\right), \label{eq:PSF}
\end{align}
where \( \o \) denotes the optical waist of our one-dimensional setup \cite{E1974, ZZOM2007}. Here we assume that our confocal volume is centered at the origin of our coordinate system.

\subsubsection{Data}\label{sec:data} 
\begin{figure}
    \label{fig:timeline}
    \includegraphics[width=6.5in]{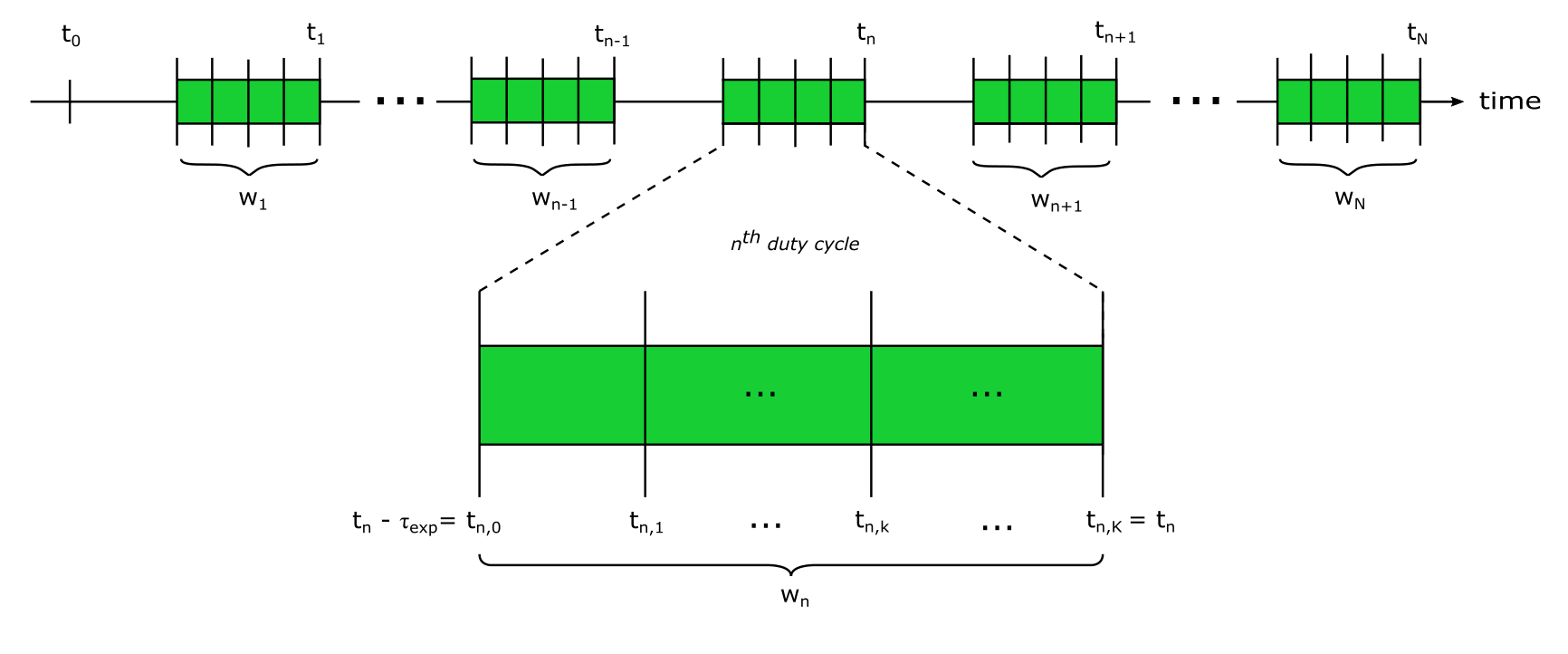}
    \center
    \caption{The measurement scale mesh is indexed by \( n=1:N \), with highlighted exposure windows and dead time periods in between. We show one such window blown-up with the superresolution subpanel mesh, indexed by \( k=1:K \). The measurement \( w_n \) is the photon count for the \( n^{\text{th}} \) duty cycle. Overall the experiment starts at \( t_0 \) and ends at \( t_N \).}
\end{figure}

The data obtained in a fluorescence spectroscopy experiment take the form of a time series \( w=w_{1:N} \). Namely, the experiment begins at time \( t_0 \) and ends at time \( t_N \) with measurements \( w_n \) indexed with \( n=1:N \) recorded at regular times \( t_n \) between them (see \cref{fig:timeline}). Each \( w_n \) corresponds to a period during which the detector collects the photons emitted from the molecule moving through the volume during the \( n^\text{th} \) duty cycle which runs from time \( t_n-\texp\) to time \( t_n \). Here \( \texp \) denotes the duration of the \emph{exposure period} during which photons are being collected, which we assume is constant across cycles. We also assume the interexposure time period, that is the dead time during which no photons are being collected, denoted \( \tdead \), is uniform over the cycles. At the end of the duty cycle we record \( w_n \), the cumulative number of photons detected during that cycle.

We describe the expected number of photons emitted during a period of active measurement by \( u_n \). This is an implicit function of \( x(t) \), i.e.~the molecule's trajectory over times \( t\in[t_n-\texp,t_n] \). We assume integrative photon detection \cite{K2021,K2021extraction} and obtain $u_n$ by
\begin{align}\label{eqn:stoch_integral}
    u_n = \int_{t_n-\texp}^{t_n} I(x(t))\;dt,
\end{align}
where \( I(x(t)) \) is the intensity of light at the point in the confocal volume corresponding to the position of the molecule at time \( t \). Although \( x(\cdot) \) is a stochastic trajectory (see below), the integral in \cref{eqn:stoch_integral} is meant in the Riemann sense \cite{S2021}.

To allow for computational implementations (see next section), for each \( n=1:N \) we consider
\begin{align}\label{eq:unifsubmesh}
    \HT_n = \{t_n-\texp=t_{n,0}<t_{n,1}<...<t_{n,k}<...<t_{n,K}=t_n\},
\end{align}
a uniform submesh with \( K \) subpanels subdividing the interval \( [t_n-\texp, t_n] \). We define a one-dimensional single molecule discrete trajectory \( q=q_{1:N,0:K} \) with initial position \( q_0 \) on the time mesh \( \HT = \cup_n \HT_n \) by
\begin{align}\label{eq:disctraj}
    q_0 &= x(t_0), & q_{n,k}&=x(t_{n,k}), & n&=1:N,& k&=0:K.
\end{align}
Due to the interexposure dead time periods, our $\HT$ is a two scale mesh characterized by durations $\tdead$ and $\texp/K$.

To simplify notation, we make the identifications \( q_0 = q_{0,K} \) and \( q_{n,K} = q_{n+1,-1} \), and define \( \tsub = \texp/K \). 
Further, we apply the {\it composite trapezoid rule} to approximate \( u_n \) by
\begin{align}
    u_n \approx \frac{\tsub}{2}\sum_{k=1}^{K}[I(q_{n,k-1})+I(q_{n,k})]. \label{eq:trap}
\end{align}
Our approximation is valid as long as \( K\gg 1 \). Accordingly, we assume \( K \) is sufficiently large and ignore any approximation error and call the resulting sequence \( u=u_{1:N} \) the {\it signal}.

Due to the shot noise contaminating photon emission and detection, fluorescence spectroscopy experiments yield measurements following Poisson statistics \cite{J2019, K2021}. So, we adopt the following sampling rule
\begin{align}
    w_n\mid q_{n,0:K}\sim\mathcal{P}(u_n), \label{eq:meas}
\end{align}
where \( \mathcal{P}(u_n) \) denotes the Poisson distribution \cite{P2023} with mean \( u_n \).

\subsubsection{Motion}
\label{sec:motion}
We model the motion of a molecule by \emph{free diffusion} \cite{J2019}. The dynamics over the time domain of interest \( [t_0,t_N] \) are governed by the formal stochastic differential equation
\begin{align}\label{eq:1}
\frac{d}{dt} x(t)&=\sqrt{2D}\;\xi(t),& t\in[t_0,t_N], 
\end{align}
where \( D \) is the molecule's diffusion coefficient. \Cref{eq:1} defines the position of the molecule at time \( t \) as a random variable
\begin{align}
    x(t)&=x(t_0)+\sqrt{2D}\int_{t_0}^t \xi(s)\;ds,& t\in[t_0,t_N],
\end{align}
where the random fluctuation \( \xi(t) \)
is implicitly defined by
\begin{align}
    \int_{t'}^{t''}\xi(t)\;dt&\sim\mathcal{N}(0,t''-t'),& t'<t'',
\end{align}
where \( \mathcal{N}(\m,v) \) is a normal distribution with mean \( \m \) and variance \( v \) \cite{J2019}.

Consequently, the discrete trajectory \( q \) of \cref{eq:disctraj} describes a Gaussian random walk, a discretization of the free diffusion model just given. We assume \( q_0=0 \), that is, the initial position is at the center of the confocal volume. Following our free diffusion model, the resulting equations that dictate the sampling rules for the position of the molecule \( q_{n,k} \) after the initial time are
\begin{align}
    q_{n,0}\mid q_{n-1,K} &\sim\mathcal{N}(q_{n-1,K}, 2D\tdead), & n=1:N, \label{eq:disctraj1}
    \\
    q_{n,k}\mid q_{n,k-1}&\sim\mathcal{N}(q_{n,k-1}, 2D\tsub), & k=1:K, \; n=1:N. \label{eq:disctraj2}
\end{align}
 We note that with the assumption \( K\gg 1 \), we assess the molecule's motion at a time resolution which is finer than that of data acquisition. Furthermore, this fine timescale means that the random vector \( q \), which represents a molecular trajectory, has highly correlated entries.

Given values for the parameters \( \Iref,\Ibg,\omega, D,\tdead,\texp\) and sizes \( N,K, \) our model can be simulated using ancestral sampling to generate synthetic data that mimic a real-life confocal fluorescence spectroscopy experiment \cite{BN2006,P2023}. The algorithms involved require direct sampling of Gaussian and Poisson variates which are commonly found in scientific software.

\subsection{Model inference}\label{sec:modinf} The probabilistic model of the previous section allows us to move towards the solution of the decoding problem \cite{P2023}. That is, given measurements \( w \), we seek to reconstruct a distribution for the underlying molecule trajectory \( q \).

With the goal of solving the decoding problem, we factorize \( P(q\mid w) \), the posterior probability density for the discretized trajectories of the molecule that are consistent with \( w \). An application of Bayes' theorem \cite{G1995} along with the model equations gives
\begin{align}\label{eq:targetpost}
    P(q\mid w)\propto P (w\mid q)P(q) = \left[\prod_{n=1}^NP(w_n\mid q_{n,1:K}) \right]\left[\prod_{n=1}^N\prod_{k=0}^KP(q_{n,k}\mid q_{n,k-1})\right],
\end{align}
where \( P (w\mid q) \) is the likelihood, determined by \cref{eq:meas}, and \( P(q) \) is the prior, determined by \cref{eq:disctraj1,eq:disctraj2}.

Note that this factorization is a product of known probability functions (i.e., Poisson and Normal), and it is proportional to the posterior with proportionality constant independent of \( q \). Since the target posterior is nonstandard in the sense that directly sampling it is not possible, this motivates the application of a Metropolis-type Markov chain Monte Carlo (MCMC) method \cite{L2001,RC1999}. Recall that the dimensionality of the space on which our posterior distribution lives is large (order \( NK \) with \( N\gg 1 \) and \( K\gg 1 \)), and that the coordinates of \( q \) are highly correlated, which results in a random walk MCMC method would converging very slowly \cite{GJM2011,B2017}. Because \( P(q\mid w) \) is smooth with respect to \( q \), we opt for HMC which is well suited for problems with such characteristics \cite{GC2011,GJM2011,B2017}.

\subsubsection{Hamiltonian Monte Carlo}\label{sec:HMC}

 We assume familiarity with MCMC  methods in which a sequence of samples \( \{q^{(j)}\}_{j=0:J} \) following prescribed statistics is obtained through a Metropolis-type accept/reject scheme. The classic reference is Hastings \cite{H1970}, and more modern treatments can be found in Liu or Robert and Casella \cite{L2001,RC1999}
 
 The generic HMC algorithm is a version of the Metropolis MCMC algorithm with a proposal rule constructed in analogy with classical Hamiltonian mechanics \cite{GC2011,GJM2011,B2017}. It leverages the gradient information of the target distribution by viewing the current sample \( q^{(j)} \) as the initial position of a virtual particle, or of an ensemble of virtual particles, governed by Hamiltonian dynamics in a virtual phase space. This analogy introduces auxiliary "momentum" variables \( \{p^{(j)}\}_{j=0:J} \) and a "total energy" Hamiltonian function with the form
\begin{align}\label{eq:energy}
    H(q,p) = V(q)+\frac{1}{2m}p\cdot p,
\end{align}
where \( V \) is a virtual "potential energy" chosen to satisfy
\begin{align}\label{eq:gibbsmeasure}
    e^{-V(q)}\propto P(q\mid w),
\end{align}
and \( m \) is a virtual "mass" which is a free parameter to be chosen. In \cref{eq:energy} the infix symbol "\( \cdot \)" denotes the Euclidean inner product. We choose a separable Hamiltonian in \cref{eq:energy} close to that in the original algorithm proposed by Duane et al. \cite{DKPR1987}, but various other formulations for the Hamiltonian have been proposed \cite{GJM2011,GC2011} and may be used as well.

We take the Hamiltonian dynamics of \( q \) and \( p \) to depend on a virtual time \( \eta \), henceforth referred to as HMC-time, as opposed to the physical time heretofore denoted by \( t \). We denote differentiation with respect to \( \eta \) using Newton's dot notation. Then, given a sample \( q^{(j)} \) in the MCMC chain and a random initial momentum \( p^{(j)}\sim \HN(0,mI) \), the proposal sample comes from considering dynamics governed by the initial value problem for Hamilton's equations
\begin{align}
    \dot q &= +H_p & q(0) &= q^{(j)} \label{eq:Ham1}\\ 
    \dot p &= -H_q & p(0) &= p^{(j)}, \label{eq:Ham2}
\end{align}
where \( H_p \) and \( H_q \) stand for the gradient of \( H \) with respect to \( p \) and \( q \), respectively. The random choice of initial momentum ensures ergodicity of our Markov chain, and the choice of the mass parameter \( m \) determines the average magnitude of the initial momentum \cite{B2017}. 

We denote by
\begin{align}
    \left(\tilde q^{(j+1)},\tilde p^{(j+1)}\right) = \phi_h^L\left( q^{(j)}, p^{(j)}\right)
\end{align}
our proposed sample. This is the terminal condition produced by numerical integration of \cref{eq:Ham1,eq:Ham2} for \( L \) successive HMC-timesteps of size \( h \) with initial conditions \( \left( q^{(j)}, p^{(j)}\right) \). For our application, numerical integration is required here, in which case a \emph{symplectic integrator} is the standard choice \cite{B2017, GJM2011, L2001}. The defining feature of such integrators is that the Jacobian \( \partial \phi_h \) of the one-step numerical integration map \( (q,p)\mapsto\phi_h(q,p) \) satisfies
\begin{align}
    \begin{bmatrix} 0 & -I \\ +I & 0 \end{bmatrix}\partial \phi_h \begin{bmatrix} 0 & +I \\ -I & 0 \end{bmatrix} &= \partial \phi_h
\end{align}
when applied to Hamiltonian systems. The \( L \)-step numerical integration map \( \phi_h^L \) defined as the \( L \)-fold composition of \( \phi_h \) inherits this property \cite{HLW2006}.

We prefer symplectic integrators because they are measure-preserving on \( (q,p) \)-space and accurate at long HMC-time scales \cite{HLW2006,LR2004}. The accuracy directly affects the acceptance ratio in the Metropolis step. After integration we accept the proposal with probability
\begin{align}
    R^{(j)} = \min\left(1,\exp(H( q^{(j)}, p^{(j)})-H(\tilde q^{(j+1)},\tilde p^{(j+1)})\right),
\end{align}
otherwise we set \( q^{(j+1)} = q^{(j)} \) and repeat starting with a new random momentum \( p^{(j+1)} \) (see \ref{alg:HMCrefl} in the supplement for our implementation).

The standard symplectic integrator used in HMC is the second-order, explicit St{\"o}rmer-Verlet numerical integrator \cite{GJM2011,BRSS2018} known as the leapfrog scheme. In this study, we denote any numerical integration maps using St{\"o}rmer-Verlet with an SV superscript. \Cref{alg:SV} in the supplement shows one step of St{\"o}rmer-Verlet integration with step size \( h \) for the Hamiltonian in \cref{eq:energy}, corresponding to the full posterior Hamiltonian system.

\subsubsection{A semi-implicit splitting method}
\label{sec:IE}

In the context of MCMC for posterior Bayesian inference, the standard choice of numerical integrator for HMC is the \( L \)-fold composition of the full posterior integrator \( \phi_h^{\text{SV}} \) implemented with the leapfrog technique. In this section we describe the main drawbacks of \( \phi_h^{\text{SV}} \) when applied to fluorescence spectroscopy data and alternatives appropriate for our setting.

Defining our Hamiltonian potential function, as in \cref{eq:gibbsmeasure}, in relation to the target posterior defined in \cref{eq:targetpost}, we get one term from the likelihood and one from the prior,
\begin{align}
    V(q) = -\log P(w\mid q)-\log P(q),
\end{align}
so we write \( V = \Vlike + \Vpri \) where
\begin{align}
    \Vlike(q) &= -\log P(w\mid q), \\
    \Vpri(q) &= -\log P(q).
\end{align}
Using \cref{eq:disctraj1,eq:disctraj2}, we determine that the potential from the prior is, up to an additive constant,
\begin{equation}\begin{split}
    \Vpri(q) &= \frac{1}{4D\tdead}\sum_{n=1}^N(q_{n,0}-q_{n-1,K})^2 +\frac{1}{4D\tsub}\sum_{n=1}^N\sum_{k=1}^K(q_{n,k}-q_{n,k-1})^2.
\end{split}\end{equation}
This potential has gradient
\begin{align}
    \grad\Vpri(q) = -\frac{1}{2D}\Delta_\t q,
\end{align}
where \( \Delta_\t \) is a \emph{tridiagonal} matrix with super- and sub-diagonals of the form
\begin{align}
    \underbrace{\frac{1}{\tdead}, \; \underbrace{\frac{1}{\tsub},\cdots}_{K \text{ times}},\cdots}_{N \text{ times}}
\end{align}
and diagonal with the structure
{\small
\begin{align}
    -\frac{1}{\tdead}, \underbrace{-\frac{1}{\tdead}-\frac{1}{\tsub}, \underbrace{-\frac{2}{\tsub},\cdots,}_{K-1 \text{ times}} -\frac{1}{\tsub}-\frac{1}{\tdead},\cdots,}_{N-1 \text{ times}} -\frac{1}{\tdead}-\frac{1}{\tsub}, \underbrace{-\frac{2}{\tsub},\cdots,}_{K-1 \text{ times}} -\frac{1}{\tsub}.
\end{align}}
Such a \( \Delta_\t \) matrix has the structure of a discrete Laplacian \cite{SW2022}. This is the discrete Laplacian corresponding to the two-scale physical time mesh \( \HT \) defined in terms of \cref{eq:unifsubmesh}. For illustrative purposes, we give an example of the two-scale discrete Laplacian \( \Delta_\tau \) for \( N=2 \) and \( K=3 \). It is convenient to introduce the notation
\begin{align}
    s_0 &= \frac{1}{\tdead},& s_1 &= \frac{1}{2}\left(\frac{1}{\tdead}+\frac{1}{\tsub}\right), & s_2 = \frac{1}{\tsub}.
\end{align}
With this notation we have
\begin{align}\label{eq:smallex}
    \Delta_\t = \left[
        \begin{array}{c|cccc|cccc}
            -s_0 &   s_0 &&&& \\ \hline
         s_0 & -2s_1 &   s_2 &&&& \\
             &   s_2 & -2s_2 &   s_2 &&&&\\
             &       &   s_2 & -2s_2 &   s_2 \\
             &       &       &   s_2 & -2s_1 &   s_0 \\ \hline
             &       &       &       &   s_0 & -2s_1 &   s_2 \\
             &       &       &       &       &   s_2 & -2s_2 &   s_2 \\
             &       &       &       &       &       &   s_2 & -2s_2 &   s_2 \\
             &       &       &       &       &       &       &   s_2 &  -s_2 \\
        \end{array}
    \right].
\end{align}

Due to the independence between the individual measurements in \cref{eq:meas}, the partial derivatives of \( \Vlike \) forming \( \grad \Vlike \) are fully decoupled. However, those of \( \Vpri \) are interrelated in a complicated way. This prompts us to treat the two terms of the full gradient \( \grad V \) \emph{differently}. Since we are sampling the posterior using Hamiltonian Monte Carlo, this motivates splitting the Hamiltonian of \cref{eq:energy} itself,
\begin{align}\label{eq:Hamsplitting}
    H(q,p) &= \Hlike(q,p) + \Hpri(q,p),
\end{align}
where \( \Hlike \) contains the terms related to the likelihood and \( \Hpri \) contains the terms related to the prior resulting from the models in \cref{sec:data} and \cref{sec:motion}, respectively.

Splitting the Hamiltonian, as in \cref{eq:Hamsplitting}, requires splitting the kinetic energy term \( T(p) = \frac{1}{2m}p\cdot p \). By analogy with Newton's second law of motion, the original Hamiltonian decomposes into the sum of a potential energy function \( V \) and a kinetic energy \( T \), so to maintain this structure in the likelihood and prior subsystems we define
\begin{align}
     \Hlike(q,p) &= \Vlike(q) + \Tlike(p), \\
    \Hpri(q,p) &= \Vpri(q) + \Tpri(p), \label{eq:priorenergy}
\end{align}
where we require that the kinetic energy splits as
\begin{align}
    \Tlike(p)+\Tpri(p) = T(p).
\end{align}
Beyond this sum structure, we have freedom to choose the decomposition of \( T \). To allow flexibility in this choice in a parsimonious way, we introduce a momentum splitting parameter \( \theta\in[0,1] \) and split the kinetic energy according to
\begin{align}
    \Tlike(p)&=\frac{1-\th}{2m}p\cdot p,&
    \Tpri(p) = \frac{\th}{2m}p\cdot p.
\end{align}

To numerically integrate the ODE in \cref{eq:Ham1,eq:Ham2}, we apply numerical integrators which leverage the split Hamiltonian structure of \cref{eq:Hamsplitting} \cite{HLW2006,LR2004}. Our method is to numerically flow the likelihood subsystem with \( \psi_h^{\text{SV}} \), a St{\"o}rmer-Verlet integration, and separately flow the prior subsystem with \( \chi_h^{\text{MP}} \), an implicit midpoint integration. We show below that the prior subsystem is stiff and so our implicit midpoint integrator for the prior subsystem is an appropriate choice.

There are different ways to compose the two flows. We choose the splitting scheme originally developed and analyzed by Strang \cite{S1968}. The numerical flow \( \phi_h^{\text{IMEX}} \) for our Hamiltonian Monte Carlo integration is a composition scheme defined by the Strang splitting
\[
    \phi_h^{\text{IMEX}} = \chi_{h/2}^{\text{MP}}\circ\psi_h^{\text{SV}}\circ\chi_{h/2}^{\text{MP}}.
\]
It is shown that the composition of symplectic maps remains symplectic, and furthermore, that the Strang splitting composition of two second order symplectic methods is again second order in \cite{HLW2006}. With this method we sidestep the problems associated with the stiffness of the discrete Laplacian, while maintaining a second order symplectic integration scheme.

To distinguish St{\"o}rmer-Verlet integration of the subsystems (integrators superscripted with "SV") from integration of the full posterior system, we introduce the abbreviation SVEX to refer to integration of the full posterior with the St{\"o}rmer-Verlet integrator \( \phi_h^{\text{SVEX}} \). In order to compare the stability of the IMEX method to the fully explicit SVEX method, it is sufficient to note that the stability regime for the SVEX method on the full posterior system is constrained by the stiffness of the prior subsystem. Thus, we can use a stability analysis of St{\"o}rmer-Verlet integration on the prior subsystem integrator, with numerical integration map \( \chi_h^{\text{SV}} \), as a surrogate for the SVEX method on the full posterior system using the \( \phi_h^{\text{SVEX}} \) integrator. We validate this surrogate approach numerically in \cref{sec:surrogate} (see below).

The relevant integrators of this study are tabulated below.

\begin{align*}
    \begin{array}{cc|c}
        \chi_h^{\text{SV}} & \text{St{\"o}rmer-Verlet} & \\
        \chi_h^{\text{MP}} & \text{Implicit midpoint} & \text{prior subsystem}\\ &&\\ \hline &&\\
        \psi_h^{\text{SV}} & \text{St{\"o}rmer-Verlet} & \text{likelihood subsystem} \\ &&\\ \hline &&\\
        \phi_h^{\text{IMEX}} & \text{Semi-implicit} & \text{full posterior system} \\
        \phi_h^{\text{SVEX}} & \text{St{\"o}rmer-Verlet}
    \end{array}
\end{align*}
Note the omission of an implicit integrator for the likelihood subsystem. We do not consider this variation, as the nonlinearities in the likelihood system induced by \cref{eq:PSF,eq:trap,eq:meas} would make their implicit integration prohibitively expensive due to the need for high-dimensional root finding at each integration step. \Cref{alg:SPLIT} of the supplement records \( L \) steps of the IMEX split HMC integration outlined above.

\subsubsection{Eigenvalue bounds for \( \Delta_\tau \)}

To justify the use of an implicit integrator for the prior subsystem requires a spectral analysis of the \( \Delta_\t \) operator, which we provide next. We recall that the two time scales appearing in the discrete trajectory \( q \) of \cref{eq:disctraj} are the larger scale of the duration between recorded measurements, indexed by a subscript \( n \), and the smaller scale of subpanels dividing the period of photon accumulation, indexed by \( k \) (see \cref{fig:timeline}). These physical time scales manifest in the structure of \( \Delta_\t \) with the appearance of \( \tdead \), the duration of interexposure dead time, and \( \tsub \), the duration of each subpanel of the collection period. We can decompose the discrete Laplacian \( \Delta_\t \) in relation to these two durations,
\begin{align}
     \Delta_\t = \frac{1}{\tdead}A + \frac{1}{\tsub}B
\end{align}
where \( A \) and \( B \) are both tridiagonal matrices of size \( (N(K+1)+1)\times(N(K+1)+1) \) (for an example see \cref{eq:smallex}).

The matrix \( A \) is a block-diagonal matrix with blocks of the form
\begin{align}
    A_1=\begin{bmatrix}
        -1 & +1 \\
         +1 &-1
    \end{bmatrix}
\end{align}
alternating with zero blocks, and \( B \) is a block diagonal matrix with one \( 0 \) for the upper-left entry and repeated blocks of the form
\begin{align}
    B_1=\begin{bmatrix}
        -1      & 1      &        \\
         1      &-2      & 1      \\
         {}     & \ddots & \ddots & \ddots \\
         {}     & {}     & 1  & -2  & 1\\
         {}     & {}     & {} &  1  & -1
    \end{bmatrix}
\end{align}
of size \( (K+1)\times(K+1) \) (refer to \cref{eq:smallex}). This explicit decomposition of \( \Delta_\t \) aids in the splitting method stability analysis to follow.

Our \( \grad \Vpri \) corresponds to action by a discrete Laplacian operator \( \Delta_\t \), which is ill-conditioned \cite{SW2022,LR2004}. If we integrate our Hamiltonian with an  \emph{explicit} symplectic integrator, like St{\"o}rmer-Verlet, the stiffness from the prior term forces the step sizes to be quite small in order to maintain stability \cite{HLW2006,BRSS2018}. We make this precise below. The naive solution would be to use an unconditionally stable symplectic method like implicit midpoint to integrate the full system, however, as in the case of implicitly integrating the likelihood subsystem alone, this would require prohibitively expensive root-finding at each integration step due to the nonlinearities in \( \Vlike \).  We address this problem with our semi-implicit splitting method justified by the following analysis on the prior term.

We denote the spectrum of a matrix \( X \) by \( \text{spec}(X) \). Due to its block diagonal structure, the spectrum of \( A \) is \( \text{spec}(A)=\{ 0\}\cup\text{spec}(A_1) \), where \( \text{spec}(A_1)=\{0,-2\} \). Likewise, \( \text{spec}(B)=\{0\}\cup\text{spec}(B_1) \). The matrix \( B_1 \) is a unit step size discrete Laplacian with Neumann boundary conditions \cite{SW2022}. A standard result in numerical spectral methods (see \cite{SW2022,LR2004}) is that the \( K+1 \) eigenvalues of \( B_1 \) are
\begin{align}
    \lambda_k &= -4\sin^2\left(\frac{\pi k}{2(K+1)}\right),& k=0:K.
\end{align}

Since \( -\Delta_\t \) is a discrete Laplacian, it is positive definite \cite{SW2022}. Recall that the maximal eigenvalue \( \l_{\text{max}} \) of a positive definite matrix \( X \) has the following variational characterization (see \cite{SW2022})
\begin{align}
    \l_{\text{max}}= \sup_{\|v\|_2=1}v^*Xv,
\end{align}
where \( \|\cdot\|_2 \) denotes the Euclidean norm.
Assuming \( K\gg 1 \), we have \( \tsub\leq\tdead \). It follows that the maximal eigenvalue \( \lambda \) of the linear map \( q\mapsto\grad \Vpri(q)=-\frac{1}{2D}\Delta_\t q \) satisfies
\begin{align}\label{eq:eigenvaluebound}
    \lambda < \frac{1}{D\tdead}+\frac{2}{D\tsub} \lesssim \frac{1}{D\tsub}
\end{align}
where the constant implicit in the \(\lesssim\) symbol is independent of all parameters.

\subsubsection{A CFL-type condition for the SVEX method}

By positive definiteness, every eigenvalue of \( \grad\Vpri \) is positive and bounded by \( \l \). This allows us to reduce the stability analysis of the prior subsystem to a simpler problem. We consider only the \emph{dominant eigenmode} of the prior subsystem corresponding to this problematic eigenvalue, simplifying the stability analysis of the St{\"o}rmer-Verlet integration of the prior subsystem to the investigation of a simple harmonic oscillator \cite{SW2022,BRSS2018}.  Under the assumption \( K\gg 1 \), the order of the eigenvalue corresponding to the dominant eigenmode is bounded by \( 1/D\tsub \), as shown in \cref{eq:eigenvaluebound}. Thus to perform this spectral analysis, it suffices to consider the simpler Hamiltonian
\begin{align}
    H_{\text{spec}}(q,p) &= -\frac{q^2}{2D\tsub}+\frac{\th p^2}{2m},& (q,p)\in \R^2,
\end{align}
where now \( q \) and \( p \) are restricted to be scalar valued functions whose dynamics are in the most problematic direction given by the eigendecomposition of \( \grad \Vpri \). This Hamiltonian gives the total energy for the undamped harmonic oscillator
\begin{align}
    \ddot q - c^2 q &= 0, &  c = \sqrt{\frac{\th }{m D \tsub}}.
\end{align}
The stability analysis of the numerical integration of such harmonic systems by St{\"o}rmer-Verlet integration has been completed \cite{BRSS2018,S2010,LR2004}, so we omit the details. For this ODE, stability of  St{\"o}rmer-Verlet with step size \( h \) requires that \( ch< 2 \). Since the value \( c \) is a function of the dominant eigenvalue of a discrete Laplacian, a CFL-type stability condition appears \cite{SW2022}. In particular, if \( h,\t,\th,m,D \) satisfy
\begin{align}\label{eq:CFL}
    \frac{h^2}{\tsub}\cdot\frac{\th}{m}\cdot\frac{1}{D}  \lesssim  1,
\end{align}
then with stepsize \( h \), St{\"o}rmer-Verlet integration of the prior subsystem is linearly stable.

\section{Results}
\label{sec:results}

All results shown below are obtained using fixed HMC parameters \( \th=1/2 \) and \( m=1 \), along with the following motion model parameters and microscope model parameters, which are typical for FCS in vitro biochemical experiments with free fluorescent dyes \cite{S2017,K2021}.

\

\begin{center}
\begin{tabular}{ |l|c|c|c| } 
 \hline
     Description & Parameter & Value & Units \\
 \hline
 diffusion coefficient & $D$ & $5\times 10^2$ & $\mu$m$^2$/sec \\
 reference photon emission rate & $\Iref$ & $5\times 10^4$ & photons/sec \\
 background photon emission rate & $\Ibg$ & $1\times 10^3$ & photons/sec \\
 optical waist & $\o$ & $2.3\times 10^{-1}$ & $\mu$m \\
 dead time duration & $\tdead$ & $1\times 10^{-6}$ & sec \\
 exposure time duration & $\texp$ & $9\times 10^{-5}$ & sec \\
 \hline
\end{tabular}
\end{center}

\

As mentioned in \cref{sec:motion}, our model allows for the generation of synthetic data provided known values for the parameters \( D,\Iref,\Ibg,\omega,
\allowbreak \tdead, \texp \) and sizes \(  N, K \). To obtain synthetic data, we first generate a molecule trajectory \( q \) by iteratively sampling the appropriate prior distributions, following \cref{eq:disctraj1,eq:disctraj2}. We then feed \( q \) into \cref{eq:trap} to generate a signal \( u = u_{1:N} \). Finally we simulate measurements \( w_{1:N} \) by sampling with parameters \( u_{1:N} \), as in \cref{eq:meas}. We use this method throughout to generate initial conditions for HMC integration. Also, in \cref{sec:infer} and \cref{sec:MCMC_efficiency} below we use such a simulation to generate the ground truth in benchmarking.

\subsection{Statistical inference}
\label{sec:infer}

Here we conduct an empirical check of statistical correctness. First we generate synthetic measurements \( w_{1:N} \), then construct a Markov chain \( q^{(0:J)} \) using \cref{alg:HMCrefl}, producing a collection of inferred trajectories given the measured photon counts. We do this for the SVEX scheme in the stable regime as well as for the IMEX scheme.

\begin{figure}
    \label{fig:statinfer}
    \includegraphics[width=6.5in]{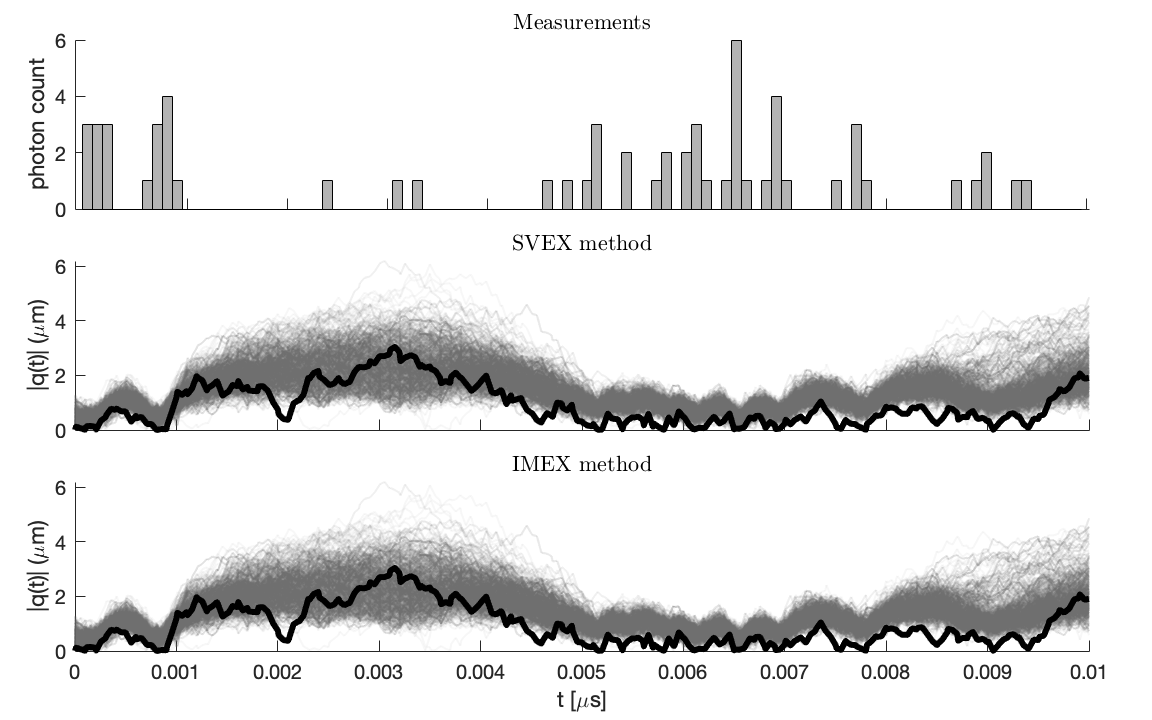}
    \center
    \caption {The top panel shows the measured photon counts. The middle and bottom panels show \( |q^{(0:J)}| \) for the SVEX scheme and the split IMEX scheme, respectively, along with the ground truth in bold.}
\end{figure}

In \cref{fig:statinfer}, we plot only the distance of the trajectory from the center of the confocal volume to better visualize the output in light of the PSF symmetry over the origin. The ground truth trajectory is plotted in bold black, while we plot 1000 MCMC trajectories for each method in gray. The two chains are both initialized from the ground truth.

The behavior is as expected, namely in the stability regime for the SVEX method, both methods produce similar output when starting from ground truth, reflecting the actual posterior statistics of \( P(q\mid w) \). We see higher variance in the sampled trajectories for both cases, for instance between 0.002 \( \m \)s and 0.004 \( \m \)s, when the ground truth trajectory is farther away from the origin. This drifting from the center of brightness shows up in the photon count as a gap in recorded photons. Note also, if the trajectory drifts close to the center of brightness, but the molecule does not emit photons, as is the case around 0.002 \( \m \)s, the inferred trajectories behave as expected, producing traces which would not likely emit photons at the 0.002 \( \m \)s mark. In other words, uninformative data produces traces with higher variance.

\subsection{Prior subsystem as stability surrogate}\label{sec:surrogate}
Here we empirically validate the use of the prior Hamiltonian subsystem as a surrogate in the stability analysis of Str{\"o}rmer-Verlet integration of the full posterior Hamiltonian system. To that end we fix the HMC integration parameters \( h \) and  \( L \), then define for a given initial position \( q_0 \) and one-step integrator \( \phi_h \) the quantity
\begin{align}
    b_{\phi_h}(q_0) = \max_{\ell=0:L}\|\phi_h^\ell(q_0)\|,
\end{align}
where \( \phi_h^0(q_0) = q_0 \). Our \( b_{\phi_h}(q_0) \) gives the uniform least upper bound for the Euclidean magnitude of a numerical trajectory starting at \( q_0 \) and integrated with \( \phi_h \) up to HMC-time \( Lh \). We expect this quantity to remain bounded for step sizes \( h \) in the stability region of the integrator.

\begin{figure}
    \label{fig:surrfig}
    \includegraphics[width=6.5in]{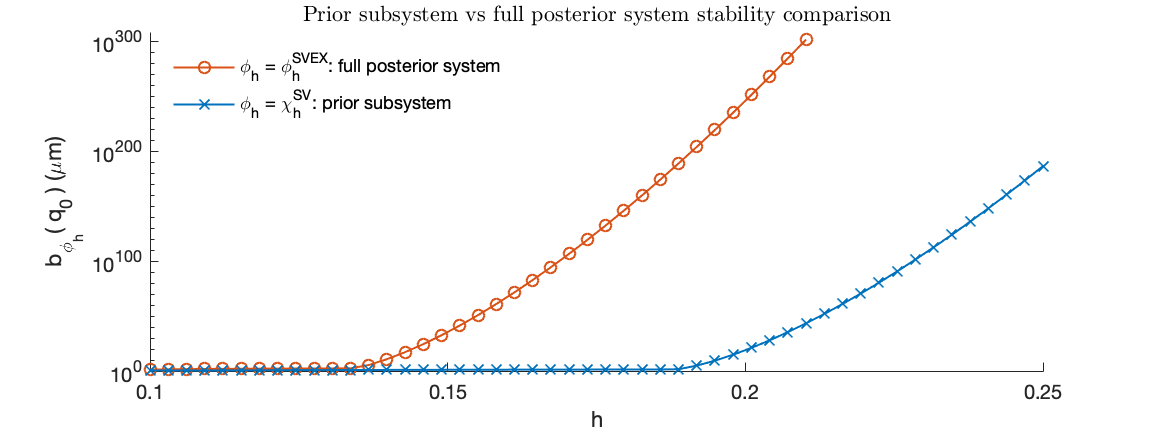}
    \center
    \caption {This figure shows \( b_{\phi_h}(q_0) \) against \( h \) on a vertical semilog scale for both the full posterior system integrated with \( \phi_h = \phi_h^{\text{SVEX}} \) (plotted in red with circular ticks), as well as the prior subsystem integrated with \( \phi_h = \chi_h^{\text{SV}} \) (plotted in blue with cross marks). The same value for \( q_0 \), was used for all trials with \( q_0 \) sampled from the prior. The measurements \( w \) used in the posterior system were generated via ancestral sampling. All other parameters were held fixed including \( N=20, K=20 \) and \( L=20 \). These plots are representative of the general behavior, with repeated simulations for randomly generated \( q_0 \) being indistinguishable at this scale.}
\end{figure}
\Cref{fig:surrfig} shows that the fully explicit integrator \( \phi_h^{\text{SVEX}} \) for the full posterior system exits its stability regime well before the fully explicit integrator \( \chi_h^{\text{SV}} \). This demonstrates that the stability analysis for the prior subsystem is an appropriate surrogate for that of the full posterior system. As expected, the terms contributed by the likelihood only increase the stiffness of the Hamiltonian system.

\subsection{Integrator stability}

In this section we numerically validate the CFL condition appearing in \cref{eq:CFL}. We test the explicit St{\"o}rmer-Verlet integrator only on the prior subsystem in order to show that the stiffness of the discrete Laplacian is a problem. In the stability regime given by the CFL condition the integrator is well behaved, in the sense that both the phase space trajectories and the total energy are bounded. However, outside of the stability regime the phase space trajectories and the total energy are unbounded. We hold all parameters and initial conditions fixed and only vary the HMC integration step size \( h \).

\begin{figure}
    \label{fig:intstab}
    \includegraphics[width=6.5in]{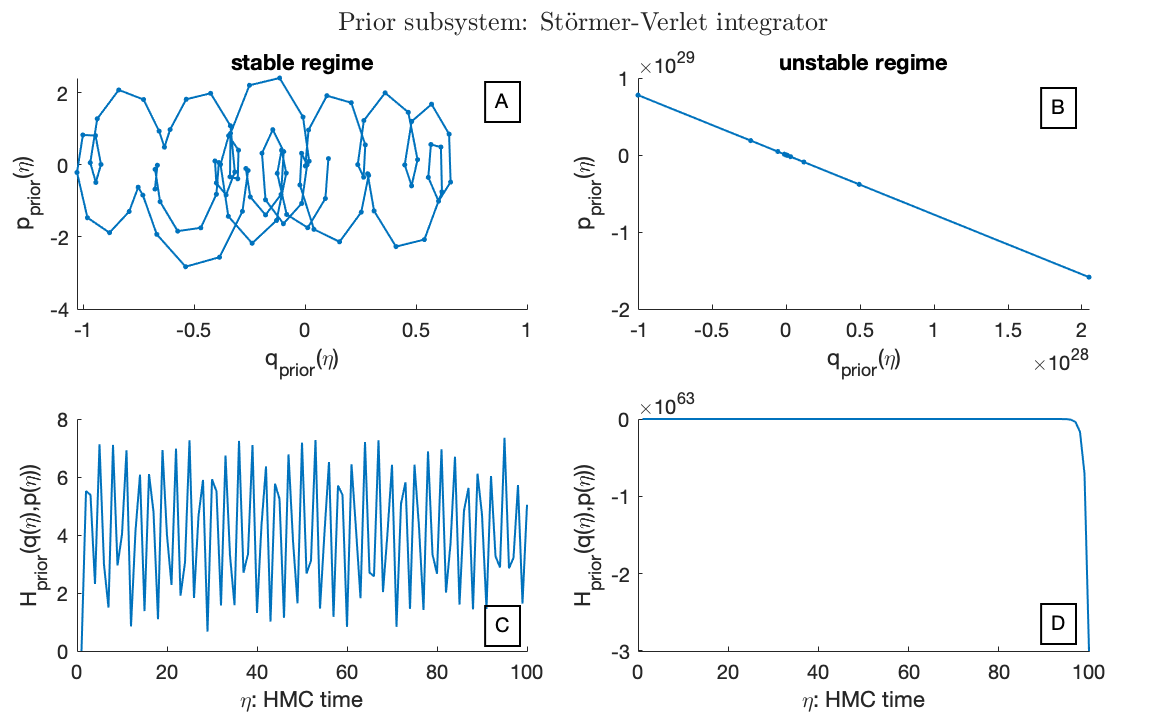}
    \center
    \caption{Panels A and B display the phase space trajectory for a single coordinate of \( (q_{\text{prior}}(\eta),p_{\text{prior}}(\eta)) \), the numerical solution to the prior subsystem integrated for \( L=100 \) steps with St{\"o}rmer-Verlet from the same initial conditions. The omitted coordinates have qualitatively similar behavior. Panels C and D show the total energy \( H_{\text{prior}} \) of the solutions sampled in panels A and B, respectively, as a function of HMC-time \( \eta \).}
\end{figure}

We remark here that, in \cref{fig:intstab}, panel A demonstrates that the phase space trajectory remains bounded with a choice of HMC step size \( h = 0.1 \) within the stability regime, while panel B shows that outside of the stability regime (\( h = 0.2 \)) the phase space trajectories are unbounded in both position and momentum. Notice the relative scale on the axes in panels A and B. Also the energy plots behave analogously. Panel C demonstrates that the total energy \( H_{\text{prior}} \) is bounded over the integration time with \( h = 0.1 \), while the total energy blows up for \( h = 0.2 \) in panel D (again, notice the relative scales of the axes).

 As demonstrated, the stiffness in the prior subsystem requires either an exceedingly small HMC step size or an integrator other than St{\"o}rmer-Verlet. Since a smaller step size translates into slower mixing of the MCMC chain given the same \( L \), we are justified in developing our alternative method. For the sake of parsimony, we choose the simplest unconditionally stable symplectic integrator, the implicit midpoint method \cite{LR2004,HLW2006}.

\subsection{MCMC efficiency}
\label{sec:MCMC_efficiency}

In this section we demonstrate how numerical instability in the prior subsystem manifests in the MCMC efficiency, that is how badly it degenerates the acceptance rate of MCMC samples. These results are concerned with sampling the full posterior, so the stability of the likelihood subsystem also effects the results. One justification for introducing the splitting scheme in the first place is that the stiffness of the prior subsystem (which incorporates the motion model in \cref{eq:disctraj1,eq:disctraj2}) and the likelihood subsystem (which incorporates the microscope and photon detector model in \cref{eq:trap,eq:meas}) may be quite different depending on the optics hardware. Our method allows for modularity in swapping out the microscope model while fixing the motion model. The following results are obtained using the confocal fluorescence spectroscopy model outlined in \cref{sec:model}.

\begin{figure}
    \label{fig:mcmceff}
    \includegraphics[width=6.5in]{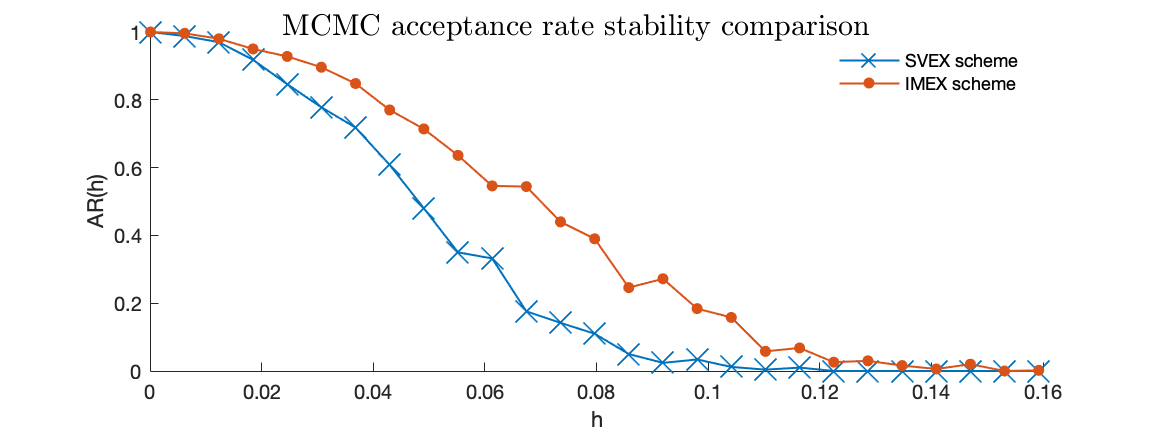}
    \center
    \caption{ Here we compute \( AR(h) \) for a range of values of \( h \) for both the SVEX numerical scheme and the IMEX splitting scheme with the implicit midpoint method integrating the prior subsystem. As expected, we see a faster decay in the MCMC stability for the explicit scheme. For example, with a step size of \( h = 0.06 \) the explicit scheme has an \( AR \) of around \( 0.38 \) while the IMEX scheme has an \( AR \) of about \( 0.63 \). }
\end{figure}

Due to the inherent stochasticity of MCMC methods we analyze multiple chains with random initial conditions generated by the model. A naive approach would be to integrate the SVEX method and the IMEX method for several choices of HMC step size \( h \) and some fixed integration time \( L \), and then compare their terminal acceptance rates, i.e., the proportion of accepted samples in relation to all proposed samples. To that end, we define \( AR_L(h) \) to be the terminal acceptance rate after \( L \) integration steps of size \( h \). It is known that the system resonances are sensitive to the total integration time \( Lh \) \cite{GJM2011,BRSS2018}, meaning that if we fix the number of HMC integration steps \( L \), there will be periodic crests in the graph of \( AR_L(h) \), even after averaging over several runs. Likewise, if the oscillations occur for some value of \( L \), they will also occur for whole number multiples of \( L \). To avoid these distortions, we define, for a fixed HMC step size \( h \), the following figure of merit for MCMC acceptance rate stability
\begin{align}
    AR(h) = \frac{1}{|\mathcal{L}|}\sum_{L\in\mathcal{L}}AR_L(h),
\end{align}
where \( \mathcal{L} \) is the set of primes less than \( 100 \). We choose the primes to avoid constructive interference of the oscillations. Fixing all other parameters while varying the HMC integration step size \( h \) yields the results displayed in \cref{fig:mcmceff}.

These results demonstrate that with our confocal fluorescence spectroscopy model, the IMEX scheme allows for larger HMC step sizes in the common range of target HMC acceptance rates \cite{BRSS2018}. Holding all else equal, larger HMC step size translates into faster MCMC convergence and better mixing time.

\subsection{Error analysis} 
\begin{figure}
    \label{fig:erran}
    \includegraphics[width=6.5in]{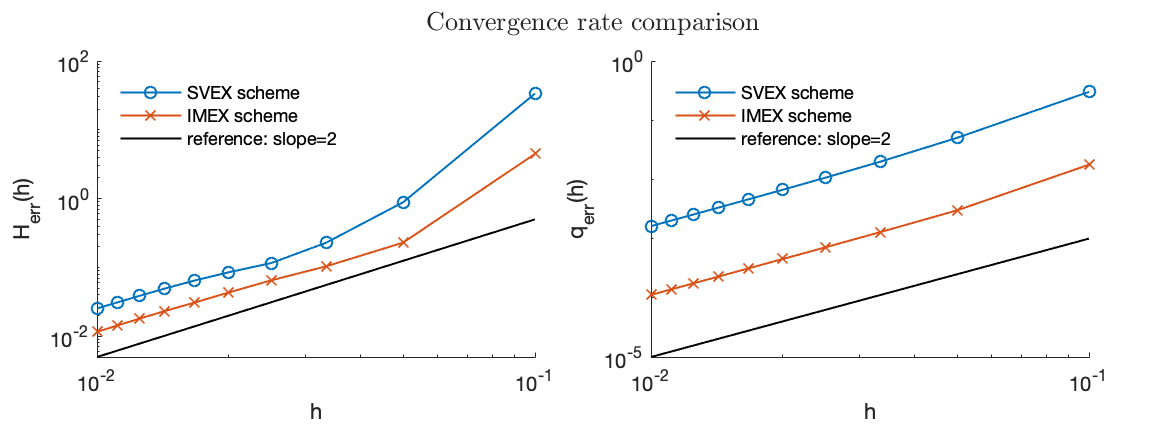}
    \center
    \caption{ The panel on the left shows the second order error convergence rate in the energy \( H_{\text{err}}(h) \) in a loglog scale.  It plots, for both schemes, the energy error \( H_{\text{err}}(h) \) against a range of \( h \) values.  A reference line of slope 2 is also plotted. The plot on the right shows the second order error convergence for the position \( q \). It plots, for both the SVEX scheme and the IMEX scheme, the terminal position error \( q_{\text{err}}(h) \) against a range of HMC step sizes \( h \) in a loglog scale. Again, a reference line of slope 2 is plotted.}
\end{figure}

As mentioned in \cref{sec:IE}, our theory predicts that both the SVEX and the IMEX scheme have second order error convergence in both the position \( q \) and the energy \( H(q,p) \) \cite{LR2004, HLW2006}. In order to validate the convergence rate we vary \( h \) while fixing \( Lh = 1 \) and integrate a fixed initial position and momentum \( (q_0,p_0) \) with a range of values of \( h \) in the stability region. We denote the terminal value of the position and momentum for a fixed \( h \) by \( (q_h,p_h) \). As a benchmark we use a step size of \( h^* = 10^{-4} \) and define the terminal position error \( q_{\text{err}}(h) \) by
\begin{align}
    q_{\text{err}}(h) = \|q_h-q_{h^*}\|_\infty.
\end{align}
The results comparing the SVEX scheme to the IMEX scheme are plotted in \cref{fig:erran}.

Since the HMC acceptance rate is computed using the terminal total energy, we are also interested in the energy error. To that end, for a fixed \( h \) we record the energy at the \( \ell^{\text{th}} \) integration step, i.e. at HMC-time \( \eta = \ell h \), and define the energy error \( H_{\text{err}}(h) \) by
\begin{align}
    H_{\text{err}}(h) = \sup_{\ell\in1:L}\Big|H(q_0,p_0)-H(q_{h}(\ell h) ,p_{h}(\ell h)\Big|
\end{align}
and plot the results for both the SVEX scheme and the IMEX scheme in \cref{fig:erran}.
As can be seen, the error in both position and energy converges quadratically.

\subsection{Complexity analysis}
\label{sec:complexity}
\begin{figure}
    \label{fig:compplot}
    \includegraphics[width=6.5in]{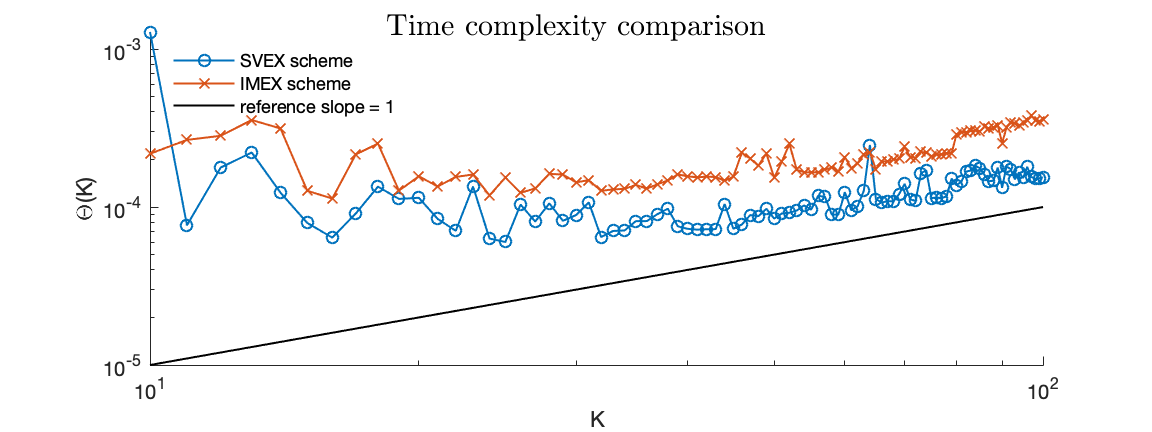}
    \center
    \caption{This plot shows, on a loglog scale, the clock time \( \Theta(K) \) against the number of subpanels \( K \) for both the SVEX (blue line with circular ticks) and the IMEX scheme (red line with cross marks) for \( K \) from 10 to 100 and fixed parameters \( L=20 \) and \( h=0.05 \). We also plot a reference line of slope 1 (solid black line). Initial conditions for the Hamiltonian system were sampled from the priors with measurements \( w \) generated via ancestral sampling, independently for each value of \( K \) but used simultaneously across the two methods for each fixed \( K \).  Both methods have linear time complexity.}
\end{figure}

Larger values of \( K \) (number of subpanels in the quadrature of \cref{eq:trap}) correspond to higher temporal resolution and model fidelity. In this section we show how computational cost is affected by \( K \). We measure wall time for the SVEX and the IMEX schemes against \( K \). To that end we consider the integration of the full Hamiltonian in both cases. The results in \cref{fig:compplot} show that the time complexity for both the SVEX and the IMEX method are comparable to leading order. Without loss, we consider a truncated molecular trajectory \( N=2 \), which does not affect the order of complexity as higher \( N \) only increases the constant multiplying the leading term.

For a fixed \( K \) we define \( \Theta(K) \) to be the wall time it takes for an \( L \)-step integration \( (L=20) \) to complete. For all simulations we used the same Apple MacBook Pro laptop computer with a 2.8 GHz Quad-Core Intel Core i7 processor. One run for each value of \( L \) and each scheme suitably demonstrates the asymptotically linear behavior for large \( K \). Both integrators are implemented using linear algebra subroutines of linear complexity, so the results are as expected \cite{SW2022}.

\section{Discussion}
\label{sec:discussion}

In this study, we introduced a high fidelity probabilistic model for scalable Bayesian data analysis of confocal fluorescence spectroscopy data, including a model for molecular motion, a model for the confocal microscope and a model for the photon accumulator. In doing so we allow for full Bayesian analysis of confocal fluorescence spectroscopy time series data, in particular the construction of a posterior for the distribution of molecule trajectories that are consistent with a given measured time series trace. An important innovation in our model is the inclusion of a superresolution parameter \( K \) which counts the number of subpanels during an exposure period and controls the discretization error or model fidelity. While only a single measurement is taken during such an exposure period (a noisy count of photons accumulated during exposure), the introduction of this subpanel parameter allows for estimates and uncertainty quantification of trajectories at time resolutions higher than the measurement time resolution which. Our estimated trajectories mimic the underlying physics of molecular motion allowing for interpretability in Chemistry and Physics applications.

In order to characterize our target posterior distribution, we constructed two scalable families of HMC sampling schemes. The IMEX split HMC algorithm, along with the CFL type condition on the stability for the SVEX scheme gives two options to practitioners. If the experimental and model parameters are in the stability regime where the CFL condition is satisfied, then the SVEX scheme may be appropriate, however the IMEX scheme allows for larger HMC step sizes while maintaining comparable MCMC acceptance rates to the SVEX scheme. Furthermore, the split structure implemented for both schemes allows practitioners to swap out microscope models without affecting the CFL stability result.

Our results demonstrate that the IMEX scheme maintains second order accuracy in HMC position and energy, while improving MCMC acceptance rates. The CFL condition predicts and our numerical experiments show that the standard fully explicit approach induces a failure mode for both sufficiently large HMC step sizes, or sufficiently many exposure window subpanels (which translates to higher resolution for predicted molecule trajectories), which motivates the use of the IMEX scheme. The stability regime for the SVEX scheme is determined by the experimental and model parameters including, among others, the molecule diffusion coefficient, as well as the HMC mass parameter. Our analysis allows the computation of a certificate for stability of the SVEX scheme in terms of these parameters which can serve as a heuristic for domain scientists and practitioners.

Both of our schemes employ the explicit St{\"o}rmer-Verlet integrator for the likelihood Hamiltonian subsystem, and since the structure of this subsystem is largely determined by the microscope model, analysis of the stiffness of the HMC subsystem induced by our particular microscope model may be a fruitful direction. Such analysis would likely give further insight into the regime of stability for both schemes and perhaps give direction for altering the microscope model to improve algorithm performance.

Possible domain specific extensions to our approach include the introduction of a 3D motion model, the modeling of multiple molecules, surface diffusion or multiple diffusion modes, and confocal geometries with Airy-type or multimodal abberated PSFs. Although a proof of principle, our work allows these extensions which may greatly improve on current FCS data analysis methods.

\bibliographystyle{plain}
\bibliography{references}

\begin{thebibliography}{10}

\bibitem{B2006}
David~Warren Ball.
\newblock {\em Field guide to spectroscopy}, volume~8.
\newblock Spie Press Bellingham, Washington, 2006.

\bibitem{B2017}
Michael Betancourt.
\newblock A conceptual introduction to {Hamiltonian Monte Carlo}.
\newblock {\em arXiv preprint arXiv:1701.02434}, 2017.

\bibitem{B1991}
EJKTDJSRL Betzig, Jay~K Trautman, TD~Harris, JS~Weiner, and RL~Kostelak.
\newblock Breaking the diffraction barrier: optical microscopy on a nanometric scale.
\newblock {\em Science}, 251(5000):1468--1470, 1991.

\bibitem{BN2006}
Christopher~M. Bishop and N.~M. Nasrabdi.
\newblock {\em Pattern Recognition and Machine Learning}.
\newblock Springer, 2006.

\bibitem{BRSS2018}
Nawaf Bou-Rabee and J.~M. Sanz-Serna.
\newblock Geometric integrators and the hamiltonian monte carlo method.
\newblock {\em Acta Numerica}, pages 1--92, 2018.

\bibitem{D2020}
Chao Du and SC~Kou.
\newblock Statistical methodology in single-molecule experiments.
\newblock {\em Statistical Science}, 35(1):75--91, 2020.

\bibitem{DKPR1987}
S.~Duane et~al.
\newblock Hybrid monte carlo.
\newblock {\em Physics Letters B}, 195:216--222, 1987.

\bibitem{D2017}
Iain~S Duff, Albert~Maurice Erisman, and John~Ker Reid.
\newblock {\em Direct methods for sparse matrices}.
\newblock Oxford University Press, 2017.

\bibitem{E1974}
Elliot~L Elson and Douglas Magde.
\newblock Fluorescence correlation spectroscopy. i. conceptual basis and theory.
\newblock {\em Biopolymers: Original Research on Biomolecules}, 13(1):1--27, 1974.

\bibitem{G1995}
Andrew Gelman, John~B Carlin, Hal~S Stern, and Donald~B Rubin.
\newblock {\em Bayesian data analysis}.
\newblock Chapman and Hall/CRC, 1995.

\bibitem{G2017}
Subhashis Ghosal and Aad~W van~der Vaart.
\newblock {\em Fundamentals of nonparametric Bayesian inference}, volume~44.
\newblock Cambridge University Press, 2017.

\bibitem{GC2011}
Mark Girolami and Ben Calderhead.
\newblock Riemann manifold langevin and hamiltonian monte carlo methods.
\newblock {\em J. R. Statist. Soc. B}, 73:123--214, 2011.

\bibitem{G2005}
Joseph~W Goodman.
\newblock {\em Introduction to Fourier optics}.
\newblock Roberts and Company publishers, 2005.

\bibitem{HLW2006}
Ernst Hairer, Christian Lubich, and Gerhard Wanner.
\newblock {\em Geometric Numerical Integration: Structure-Preserving Algorithms for Ordinary Differential Equations}.
\newblock Springer, second edition, 2006.

\bibitem{H1970}
W.~K. Hastings.
\newblock Monte carlo sampling methods using markov chains and their applications.
\newblock {\em Biometrika}, 57(1), 1970.

\bibitem{J2019}
S.~Jazani et~al.
\newblock An alternative framework for fluorescence correlation spectroscopy.
\newblock {\em Biophysical Journal}, 116, 2019.

\bibitem{J2019method}
Sina Jazani, Ioannis Sgouralis, and Steve Press{\'e}.
\newblock A method for single molecule tracking using a conventional single-focus confocal setup.
\newblock {\em The Journal of chemical physics}, 150(11), 2019.

\bibitem{J2022}
Sina Jazani, Lance~WQ Xu, Ioannis Sgouralis, Douglas~P Shepherd, and Steve Press{\'e}.
\newblock Computational proposal for tracking multiple molecules in a multifocus confocal setup.
\newblock {\em ACS photonics}, 9(7):2489--2498, 2022.

\bibitem{K2021extraction}
Zeliha Kilic, Ioannis Sgouralis, Wooseok Heo, Kunihiko Ishii, Tahei Tahara, and Steve Press{\'e}.
\newblock Extraction of rapid kinetics from smfret measurements using integrative detectors.
\newblock {\em Cell Reports Physical Science}, 2(5), 2021.

\bibitem{K2021}
Zeliha Kilic, Ioannis Sgouralis, and Steve Press{\'e}.
\newblock Generalizing hmms to continuous time for fast kinetics: Hidden markov jump processes.
\newblock {\em Biophysical journal}, 120(3):409--423, 2021.

\bibitem{K2020}
David~R Klein.
\newblock {\em Organic chemistry}.
\newblock John Wiley \& Sons, 2020.

\bibitem{L2006}
Joseph~R Lakowicz.
\newblock {\em Principles of fluorescence spectroscopy}.
\newblock Springer, 2006.

\bibitem{L2007}
Antony Lee, Konstantinos Tsekouras, Christopher Calderon, Carlos Bustamante, and Steve Press{\'e}.
\newblock Unraveling the thousand word picture: an introduction to super-resolution data analysis.
\newblock {\em Chemical reviews}, 117(11):7276--7330, 2017.

\bibitem{LR2004}
Benedict Leimkuhler and Sebastian Reich.
\newblock {\em Simulating Hamiltonian Dynamics}.
\newblock Cambridge University Press, 2004.

\bibitem{L2001}
J.~S. Liu.
\newblock {\em Monte Carlo Strategies in Scientific Computing}.
\newblock Springer, 2001.

\bibitem{L2015}
Zhe Liu, Luke~D Lavis, and Eric Betzig.
\newblock Imaging live-cell dynamics and structure at the single-molecule level.
\newblock {\em Molecular cell}, 58(4):644--659, 2015.

\bibitem{M1974}
Douglas Magde, Elliot~L Elson, and Watt~W Webb.
\newblock Fluorescence correlation spectroscopy. ii. an experimental realization.
\newblock {\em Biopolymers: Original Research on Biomolecules}, 13(1):29--61, 1974.

\bibitem{M1997}
Donald~Allan McQuarrie and John~Douglas Simon.
\newblock {\em Physical chemistry: a molecular approach}, volume~1.
\newblock University science books Sausalito, CA, 1997.

\bibitem{GJM2011}
Radford Neal.
\newblock {\em MCMC using Hamiltonian dynamic, Ch. 5 of Handbook of Markov Chain Monte Carlo}.
\newblock CRC Press, 2011.

\bibitem{P2023}
Steve Press{\'e} and Ioannis Sgouralis.
\newblock {\em Data modeling for the sciences: applications, basics, computations}.
\newblock Cambridge University Press, 2023.

\bibitem{R1993}
R~Rigler, {\"U}lo Mets, J~Widengren, and P~Kask.
\newblock Fluorescence correlation spectroscopy with high count rate and low background: analysis of translational diffusion.
\newblock {\em European Biophysics Journal}, 22:169--175, 1993.

\bibitem{RC1999}
Christian~P. Robert and G.~Casella.
\newblock {\em Monte Carlo Statistical Methods}.
\newblock Springer, 1999.

\bibitem{S2023III}
Matthew Safar, Ayush Saurabh, Bidyut Sarkar, Mohamadreza Fazel, Kunihiko Ishii, Tahei Tahara, Ioannis Sgouralis, and Steve Press{\'e}.
\newblock Single-photon smfret. iii. application to pulsed illumination.
\newblock {\em Biophysical Reports}, 2(4), 2023.

\bibitem{S2019}
Bahaa~EA Saleh and Malvin~Carl Teich.
\newblock {\em Fundamentals of photonics}.
\newblock john Wiley \& sons, 2019.

\bibitem{SW2022}
Abner~J. Salgado and Steven~M. Wise.
\newblock {\em Classical Numerical Analysis}.
\newblock Cambridge University Press, 2022.

\bibitem{S2023I}
Ayush Saurabh, Mohamadreza Fazel, Matthew Safar, Ioannis Sgouralis, and Steve Press{\'e}.
\newblock Single-photon smfret. i: Theory and conceptual basis.
\newblock {\em Biophysical Reports}, 3(1), 2023.

\bibitem{S2023II}
Ayush Saurabh, Matthew Safar, Mohamadreza Fazel, Ioannis Sgouralis, and Steve Press{\'e}.
\newblock Single-photon smfret: Ii. application to continuous illumination.
\newblock {\em Biophysical Reports}, 3(1), 2023.

\bibitem{S2021}
Ren{\'e}~L Schilling.
\newblock {\em Brownian motion: a guide to random processes and stochastic calculus}.
\newblock Walter de Gruyter GmbH \& Co KG, 2021.

\bibitem{S2010}
Tamar Schlick.
\newblock {\em Molecular modeling and simulation: an interdisciplinary guide}, volume~2.
\newblock Springer, 2010.

\bibitem{S2018}
Ioannis Sgouralis, Shreya Madaan, Franky Djutanta, Rachael Kha, Rizal~F Hariadi, and Steve Press{\'e}.
\newblock A bayesian nonparametric approach to single molecule forster resonance energy transfer.
\newblock {\em The Journal of Physical Chemistry B}, 123(3):675--688, 2018.

\bibitem{S2023}
Ioannis Sgouralis, Lance~WQ Xu, Ameya~P Jalihal, Nils~G Walter, and Steve Press{\'e}.
\newblock Bnp-track: A framework for superresolved tracking.
\newblock {\em bioRxiv}, 2023.

\bibitem{S2017}
Hao Shen, Lawrence~J Tauzin, Rashad Baiyasi, Wenxiao Wang, Nicholas Moringo, Bo~Shuang, and Christy~F Landes.
\newblock Single particle tracking: from theory to biophysical applications.
\newblock {\em Chemical reviews}, 117(11):7331--7376, 2017.

\bibitem{S1968}
G.~Strang.
\newblock On the construction and comparison of difference schemes.
\newblock {\em SIAM Journal of Numerical Analysis}, 5:506--517, 1968.

\bibitem{T2020}
Meysam Tavakoli, Sina Jazani, Ioannis Sgouralis, Omer~M Shafraz, Sanjeevi Sivasankar, Bryan Donaphon, Marcia Levitus, and Steve Press{\'e}.
\newblock Pitching single-focus confocal data analysis one photon at a time with bayesian nonparametrics.
\newblock {\em Physical review X}, 10(1):011021, 2020.

\bibitem{W1996}
R.~H. Webb.
\newblock Confocal optical microscopy.
\newblock {\em Reports on Progress in Physics}, 59, 1996.

\bibitem{Y2011}
Robert~J Young and Peter~A Lovell.
\newblock {\em Introduction to polymers}.
\newblock CRC press, 2011.

\bibitem{ZZOM2007}
B.~Zhang, J.~Zerubia, and J.~Olivo-Marin.
\newblock Gaussian approximations of fluorescence microscope point-spread function models.
\newblock {\em Applied Optics}, 46(10), 2007.

\end{thebibliography}

\newpage

\setcounter{section}{0}
\renewcommand{\thesection}{\Alph{section}}
\section{Supplement}

This supplement consists of implementation details of the algorithms described in the body of the present study.

\subsection{HMC/reflection-within-Gibbs}

The HMC update for \( q^{(j)} \) described in \cref{sec:HMC} is implemented as the first step in a two-step Gibbs sampling scheme. The model assumption of one dimensional diffusion introduces a bimodality in the target posterior. Namely, due to the PSF symmetry in \cref{eq:PSF} any segment of a molecule trajectory on one side of the center of the confocal volume is expected to produce measurements identical to the same trajectory segment reflected across the center. This bimodality can trap the sampler in cases where \( q \) crosses the the center of the confocal volume. For this reason we include the following second step in the Gibbs sampling structure. Given \( q^{(j)} \), with equal probability for two cases, a Metropolis update is performed. In case 1 (resp. case 2), for \( n=1:N,\,k=1:K \) propose \( \hat r_{n,k}(q^{(j)}) \) (resp. \( \check r_{n,k}(q^{(j)}) \)) defined to be \( q^{(j)} \) with a reflection of the head of the trajectory up to time \( t_{n,k} \) (resp. of the tail after time \( t_{n,k} \)) across the center of the confocal volume, and either accept or reject according to a Metropolis scheme. For consistency of notation, we will assume this reflection update happens "in place," with the output of the reflection update for the sample \( q^{(j)} \) also being denoted \( q^{(j)} \). We record the HMC/reflection-within-Gibbs algorithm here.
\begin{algorithm}
\caption{HMC/reflection-within-Gibbs}\label{alg:HMCrefl}
\begin{algorithmic}
\STATE \textbf{Input:} $q^{(0)}$ (Required to be in the support of the target distribution)
\STATE \textbf{Output:} $\{q^{(j)}\}_{j=1:J}$
\FOR{$j=0:J-1$}
    \STATE ${}$
    \STATE $\text{\# HMC update}$
    \STATE $p^{(j)}\sim\HN(0,mI)$
    \STATE $(\tilde q^{(j+1)},\tilde p^{(j+1)}) \gets \phi_h^L(q^{(j)},p^{(j)})$
    \STATE $\gamma^{(j)} \sim \text{Bernoulli}(R^{(j)})$
    \STATE $q^{(j+1)} \gets \gamma^{(j)} \tilde q^{(j+1)} + (1-\gamma^{(j)}) q^{(j)}$
    \STATE ${}$
    \STATE $\text{\# Reflection update}$
    \STATE $\alpha^{(j+1)}\sim\text{Bernoulli}(1/2)$
    \FOR{$n=1:N,\,k=1:K$}
        \STATE $r_{n,k} = \alpha^{(j+1)}\hat r_{n,k}+(1-\alpha^{(j+1)})\check r_{n,k}$
        \STATE $\beta^{(j+1)}\sim\text{Bernoulli}(\min(1,P(r_{n,k}(q^{(j+1)}))/P(q^{(j+1)})))$
        \STATE $q^{(j+1)}\gets \beta^{(j+1)} r_{n,k}(q^{(j+1)})+(1-\beta^{(j+1)})q^{(j+1)}$
    \ENDFOR
\STATE
\ENDFOR
\end{algorithmic}
\end{algorithm}

\subsection{Leapfrog integration of a separable Hamiltonian}
\begin{algorithm}
\caption{St{\"o}rmer-Verlet numerical integrator $\phi_h^{\text{SV}}$\\(one step, step size $h$, full posterior system)}\label{alg:SV}
\begin{algorithmic}
\STATE \textbf{Input:} $(q^0,p^0)$
\STATE \textbf{Output:} $(q^1,p^1) = \phi_h^{\text{SV}}(q^0,p^0)$
\STATE $p^{1/2} \gets p^0-\frac{h}{2}\grad V(q^0)$
\STATE $q^{1} \gets q^0+\frac{h}{m}p^{1/2}$
\STATE $p^1 \gets p^{1/2}-\frac{h}{2}\grad V(q^1)$
\STATE 
\end{algorithmic}
\end{algorithm}

For \( L \)-step integration, using the so-called leapfrog technique gives a significant computational savings \cite{HLW2006,GJM2011,BRSS2018}. This modification consists of altering the inner iterations, folding the half-step momentum update for \( p^{\ell} \) into the half-step momentum update for \( p^{\ell+1/2} \),
\begin{align}
    p^{\ell+1/2} &\leftarrow p^{\ell-1/2}-h\grad V(q^\ell), & \ell=2:L-1,
\end{align}
thus eliminating an explicit update for \( p^{\ell},\;\ell=2:L-1 \).

\subsection{SVEX and IMEX split integration schemes}

\begin{algorithm}
\caption{IMEX Symplectic split integrator (L steps, full posterior system)}\label{alg:SPLIT}
\begin{algorithmic}
\STATE \textbf{Input:} $(q^0,p^0)$
\STATE \textbf{Output:} $(q^L,p^L)=(\phi_h^{\text{IMEX}})^L(q^0,p^0)$
\STATE $(q^{1/4},p^{1/4})\gets \chi_{h/2}^{\text{MP}}(q^0,p^0)$
\STATE $(q^{3/4},p^{3/4})\gets \psi_h^{\text{SVEX}}(q^{1/4},p^{1/4})$
\FOR{$\ell=1:L-1$}
    \STATE $(q^{\ell+1/4},p^{\ell+1/4})\gets \chi_h^{\text{MP}}(q^{\ell-1/4},p^{\ell-1/4})$
    \STATE $(q^{\ell+3/4},p^{\ell+3/4})\gets \psi_h^{\text{SVEX}}(q^{\ell+1/4},p^{\ell+1/4})$
\ENDFOR
\STATE $(q^L,p^L)\gets \chi_{h/2}^{\text{MP}}(q^{L-1/4},p^{L-1/4})$
\end{algorithmic}
\end{algorithm}

In \cref{alg:SPLIT} we take advantage of the split structure of \cref{eq:Hamsplitting}. In particular, we use the fact that the operator \( \phi_h^{\text{IMEX}} \) telescopes similar to the leapfrog technique, i.e.,
\begin{align}
    (\phi_h^{\text{IMEX}})^L = \chi_{h/2}^{\text{MP}}\circ(\psi_h^{\text{SVEX}}\circ\chi_h^{\text{MP}})^{L-1}\circ\psi_h^{\text{SVEX}}\circ\chi_{h/2}^{\text{MP}}.
\end{align}
The explicit form of the \( \psi_h^{\text{SVEX}} \) subroutine for explicit integration of the likelihood subsystem is given in \cref{alg:SVlike}.
\begin{algorithm}
\caption{St{\"o}rmer-Verlet numerical integrator $\psi_h^{\text{SVEX}}$\\(one step, step size $h$, likelihood subsystem)}\label{alg:SVlike}
\begin{algorithmic}
\STATE \textbf{Input:} $(q^0,p^0)$
\STATE \textbf{Output:} $(q^1,p^1)=\psi_h^{\text{SVEX}}$
\STATE $p^{1/2}\gets p^0-\frac{h}{2}\grad \Vlike (q^0)$ \\
\STATE $q^1 \gets q^0 + h \grad \Tlike(p^{1/2})$ \\
\STATE $p^1 \gets p^{1/2} - \frac{h}{2}\grad\Vlike(q^1)$
\end{algorithmic}
\end{algorithm}

Next we describe our method for advancing one step of \( \chi_{h/2}^{\text{MP}} \) the implicit midpoint integration of the prior subsystem. By substituting the Hamiltonian for the prior subsystem \cref{eq:priorenergy} into Hamilton's equations \cref{eq:Ham1,eq:Ham2}, we obtain the ODE system
\begin{align}
    \dot q &= \frac{\th}{m}p, \\
    \dot p &= \frac{1}{2D}\Delta_{\t}q,
\end{align}
which has, for an implicit midpoint integration scheme, the update rule
\begin{align}
    q^1&=q^0+\frac{h}{2}\left(\frac{\th}{m}p^1+\frac{\th}{m}p^0\right), \\
    p^1&=p^0+\frac{h}{2}\left(\frac{1}{2D}\Delta_{\t}q^1+\frac{1}{2D}\Delta_{\t}q^0\right).
\end{align}
This reduces to the solution of two \emph{tridiagonal} linear systems,
\begin{align}
    \left(I-\frac{\th h^2}{8D m}\Delta_{\t}\right)q^1 &= \left(I+\frac{\th h^2}{8D m}\Delta_{\t}\right)q^0+\frac{\th h}{m}p^0,\\
    \left(I-\frac{\th h^2}{8D m}\Delta_{\t}\right)p^1 &= \left(I+\frac{\th h^2}{8D m}\Delta_{\t}\right)p^0+\frac{h}{2D}\Delta_{\t}q^0,
\end{align}
each of whose right hand side involves tridiagonal matrix-vector multiplication.
Our method uses the Thomas algorithm for tridiagonal systems to solve for \( (q^1,p^1) \), and the tridiagonal matrix-vector multiplication is implemented with the standard sparse method \cite{D2017, SW2022}.

\end{document}